\documentclass[12pt, letterpaper]{article}
\usepackage{fullpage}
\usepackage[top=0.25cm, left=0.8in, right=0.8in, bottom=1.3in]{geometry}
\usepackage{amsmath,amsthm,amsfonts,amssymb,amscd}
\usepackage{lastpage}
\usepackage{enumerate}

\usepackage{fancyhdr}
\usepackage{mathrsfs}
\usepackage[dvipsnames]{xcolor}
\usepackage{listings}
\usepackage{tikz}
\usetikzlibrary{decorations.markings}
\usepackage{mathtools}
\usepackage{breqn}
\usepackage{tensor}
\usepackage{parskip}
\usepackage{subcaption}
\usepackage{tikz}
\usepackage{float} \usepackage{chngcntr}
\usepackage[allcolors = blue,colorlinks]{hyperref}
\counterwithin{figure}{section}

\theoremstyle{plain}
\newtheorem{theorem}{Theorem}[section]

\newtheorem{cor}[theorem]{Corollary}
\newtheorem{prop}[theorem]{Proposition}
\newtheorem{lemma}[theorem]{Lemma}

\theoremstyle{definition}
\newtheorem{definition}[theorem]{Definition}

\newtheorem{example}[theorem]{Example}
\newtheorem{remark}[theorem]{Remark}

\newtheorem*{convention*}{Convention}
\begin{document}
\title{On Medial Quandle and Detecting Causality}
\author{Hongxu Chen \date{} } \maketitle

\begin{abstract}
We investigate the capability of medial quandles to detect causality in (2+1) dimensional globally
hyperbolic spacetime by determining if their coloring link invariants can distinguish between the connected sum of two
Hopf links and an infinite series of relevant three-component links constructed by Allen and
Swenberg in 2020, who suggested that any link invariant must be able to distinguish those
links for them to detect causality in the given setting. We use algebraic methods to obtain a comprehensive analysis of the ability of medial quandles to complete this task. In particular, we obtain firmly useful rules to determine whether a given medial quandle is capable of distinguishing the links, and find a method of constructing a large family of medial quandle that are able to distinguish the connected sum of two Hopf links and every link in the Allen-Swenberg infinite series , which therefore gives powerful candidates for detecting causality in the spacetime of our setting.
\end{abstract}

\tableofcontents

\section{Introduction} 
Mathematical knots are embeddings of circles in 3-dimensional Euclidean space. Two knots are equivalent if and only if their diagrams can be transformed to each other via a series of Reidemeister moves and planar isotopy. A link is a collection of knots that are intertwined together. A knot or link invariant is a function that assigns the same value to equivalent knots/links, which require it to be preserved under all Reidemeister moves. Using link invariant to distinguish link has been shown to relate to detection of \textit{causality} in spacetime with the following work from former researchers: \par
Let $X$ be a $(2+1)$-dimensional globally hyperbolic spacetime with Cauchy surface $\sum$ homeomorphic to $\mathbb{R}^2$. Chernov and Nemirovski \cite{CN} proved the Low conjecture \cite{L}, which states that as long as $\sum$ is not homeomorphic to $S^2$ or $\mathbb{R}P^2$, then two events $x,y\in X$ are causally related if and only if their \textit{skies} $S_x\cup S_y$ are topologically linked. Two causally unrelated points in spacetime give a link isotopic to the connected sum of two Hopf link, denoted as $L_{2H}$ in this paper. This result gave rise to the possibility that the causality of two events in their skies can be detected by link invariants.\par
 Allen and Swenberg \cite{AS} investigated a natural question: whether $S_x$ and $S_y$, and thus the causality, can be detected by link invariant. If there exist a three-component link consisted of one unknot and other two individually deformable to the interior of the longitudes of the solid torus, $S^1\times R^2$, whose link invariant is the same as $L_{2H}$, then it is said that the link invariant can not completely determine causality in $X$. They showed that Conway polynomial cannot detect causality in that spacetime, as there exist an infinite family of three-components two-skies link (called Allen-Swenberg links in this paper), which are likely skies of non-causally related events, that can not be distinguished from  the connected sum of two Hopf links by Conway polynomial.\par

\textit{Quandle} is a type of algebraic structure. \textit{Medial quandle} is a specific type of quandle. In this paper, we use \textit{medial quandle coloring polynomial} as our link invariant to study the question of causality detection with the Allen-Swenberg links.
 
We organize this paper as follows: Section \ref{section:2} presents some background of causality in spacetime and its connection with topology, and  introduces  link invariant generated by \textit{quandle}.
In Section \ref{section 3}, we present
a comprehensive analysis of the ability of medial quandles to distinguish between $L_{2H}$ and the
Allen-Swenberg series of links. We use the whole Section \ref{sec4} to prove propositions stated in Section \ref{section 3}. Section \ref{secconclusion} summarizes our main results.

\section{Background on Causality in Spacetime and Quandle Coloring Link Invariants }\label{section:2}

\begin{definition}
Let \((X,g)\) be a time-oriented Lorentzian manifold of dimension \(m+1\). A Cauchy surface is a closed achronal set \(\Sigma \subset X\) such that every inextendible causal curve in \(X\) intersects \(\Sigma\) exactly once. A spacetime that admits a Cauchy surface is called globally hyperbolic.

 \end{definition}

 \begin{definition}
     Let \(N\) denote the set of all unparameterized, future-directed null geodesics in \(X\). This space \(N\) has a natural structure of a contact manifold. In fact, \(N\) can be identified with the spherical cotangent bundle \(ST^*\Sigma\) of a Cauchy surface \(\Sigma\).

For an event \(x \in X\), the sky \(S_x\) is the subset of \(N\) consisting of all light rays passing through \(x\):
\[
S_x = \{\gamma \in N \mid x \in \gamma\}.
\]
The sky \(S_x\) is an embedded Legendrian submanifold of \(N\), diffeomorphic to the sphere \(S^{m-1}\). In the case \(m=2\), this is a circle \(S^1\).
 \end{definition}

 \begin{theorem}[Low's Conjecture; proved by Chernov--Nemirovski \cite{CN}]
Let \(X\) be a \((2+1)\)-dimensional globally hyperbolic spacetime with Cauchy surface \(\Sigma\) whose universal cover is not diffeomorphic to \(S^2\) or \(\mathbb{R}P^2\). Two events \(x, y \in X\) are causally related if and only if their skies \(S_x\) and \(S_y\) are topologically linked in \(N\).
\end{theorem}

\begin{theorem}[Legendrian Low's Conjecture; formulated by Nat\'ario--Tod \cite{NT}, proved by Chernov--Nemirovski \cite{CN}]
Let \(X\) be a globally hyperbolic spacetime with Cauchy surface diffeomorphic to an open subset of \(\mathbb{R}^m\). Two events \(x, y \in X\) are causally related if and only if the Legendrian link \(S_x \sqcup S_y\) is non-trivial in the contact manifold \(N\).
\end{theorem}

These results provide the geometric and topological foundation for using link invariants to detect causality in globally hyperbolic spacetimes.

 \begin{definition}

 A quandle $Q(*)$ is a set $Q$ equipped with a binary operation $*: (x,y)\to x*y$ such that the operation satisfy the following axioms, $\forall x, y, z\in Q:$
  \begin{enumerate}
    \item  $x*x=x$, 
    \item  There exists a unique $q\in Q$ such that $q*y=x$. We can then define the inverse operation $*^{-1}$ as $x*^{-1}y=q$, and 
    \item  $(x*y)*z=(x*z)*(y*z)$. \\ 
\end{enumerate}
\end{definition}

\begin{definition}
    A quandle $Q(*)$ is called medial if $(x_1*y_1)*(x_2*y_2)=(x_1*x_2)*(y_1*y_2), \forall x_1,y_1,x_2,y_2\in Q$. \\ \end{definition} 

    \begin{example} \label{DofAlexan}
    An Alexander quandle is a quande constructed over any Abelian group $(G,+)$ along a group automorphism $f$, where the quandle operation is defined by $a*b=(1-f)(a)+f(b)$. It is always a medial quandle because $\forall a_1,b_1,a_2,b_2\in Q$,
    \begin{equation*}
\begin{split}
(a_1*b_1)*(a_2*b_2) & = ((1-f)(a_1)+f(b_1))*((1-f)(a_2)+f(b_2)) \\
 & = ((1-f)^2(a_1)+(f-f^2)(b_1))+((f-f^2)(a_2)+f^2(b_2)) \\ 
 & =((1-f)^2(a_1)+(f-f^2)(a_2))+((f-f^2)(b_1)+f^2(b_2)) \\
 & = (a_1*a_2)*(b_1*b_2).
\end{split}
\end{equation*}
    
    %\\ {\centering $(a_1*b_1)*(a_2*b_2)=((1-f)(a_1)+f(b_1))*((1-f)(a_2)+f(b_2))=((1-f)^2(a_1)+(f-f^2)(b_1))+((f-f^2)(a_2)+f^2(b_2))=((1-f)^2(a_1)+(f-f^2)(a_2))+((f-f^2)(b_1)+f^2(b_2))=(a_1*a_2)*(b_1*b_2)$.} 
    \end{example}

\begin{definition} \label{def_fijquandle}
Let $\{X_i\}_{i \in \Lambda}$ be a family of non-empty sets, $X = \bigsqcup_{i \in \Lambda} X_i$ , and let $f_{i,j} : X_i \to X_i$ for each $i, j \in \Lambda$. Define the binary operation $* : X \times X \to X$ as follows: for any $x_i \in X_i$ and $x_j \in X_j$, \[x_i * x_j = f_{i,j}(x_i).\]\\
 $(X,*)$ is a quandle if it satisfies the following conditions, for all $i,j,k\in \Lambda:$
\begin{enumerate}
    \item $f_{i,i}=id$
    \item $f_{i,j}$ are bijective
    \item $f_{i,j}f_{i,k}=f_{i,k}f_{i,j}$ \end{enumerate}
This type of quandle is called $f_{ij}$ quandle, introduced in \cite{MP}.

\end{definition}

\begin{definition}
    Given two quandles $A(*)$ and $B(\times)$, a homomorphism from $A(*)$ to $B(\times)$, $f\in Hom(A(*),B(\times))$, is a function  $f:A\to B$ such that $f(x*y)=f(x)\times f(y)~\forall x,y\in A$. \\
    \end{definition}

\begin{definition}
 For an oriented knot or link diagram with arc-set $L$, its fundamental quandle $(L,*)$ is the quandle generated by $L$ subjected to the system of relations formed  by crossing relations given by each crossing in the link diagram using the convention shown in Figure \ref{crossrelation}: \\
 \begin{figure}[h]
     \centering
     \includegraphics[width=0.29\linewidth]{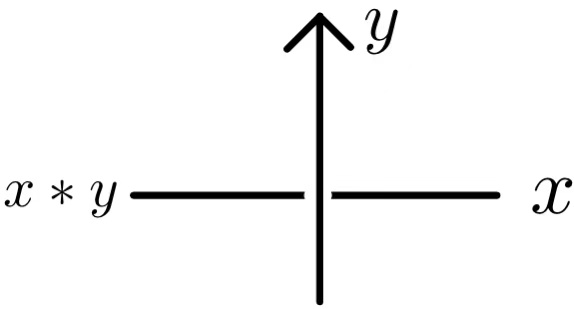}
     \caption{Crossing relation for quandle colorings. }
     \label{crossrelation}
 \end{figure} \end{definition}

\begin{theorem} \cite{J}
The fundamental quandle is a complete link invariant.
\end{theorem} 

\begin{definition} \label{countinginva}
     Given an oriented link $L$ with its fundamental quandle $Q(L,*)$ and a finite quandle $X(\times)$, the coloring space is $Hom(Q(L,*),X(\times))$, the set of homomorphisms from $Q(L,*)$ to $X(\times)$. We say that an element in the space is a coloring on $L$ by quandle $X$. Each coloring is a way of assigning an element in $X$ to each arc in $L$ such that these elements satisfy the crossing relations in $L$ but with the quandle operation in $X(\times)$. The \textit{quandle coloring counting invariant} is $|Hom(Q(L,*),X(\times))|$, the cardinality of the coloring space. \\ \end{definition} 
     \begin{remark}
Given a link $L$ and a quandle $Q$, we often represent a coloring graphically by labeling each arc in the link diagram with the element in $Q$ that is the image of the arc under the coloring.
\\The coloring space from a finite quandle $X$ on any links or knots $Q(L)$ always contains the set $C_T=\{f:Q(L)\to X,~a\to x~|~x\in X \}$, namely the set of constant functions from $Q(L)$ to $x$, where $|C_T|=|X|$. An element in $C_T$ is called the trivial coloring, sending every arc in the link or knot to the same element in the quandle. \\

\end{remark}

\begin{definition} \label{colorpolynomial}
    Introduced by Nelson \cite{NN}, using the same setting as in Definition \ref{countinginva}, the enhanced coloring invariant is a polynomial in $q$ assigned to the coloring space: \[\Phi(L,X)=\sum_{f\in Hom(Q(L),X(\times))} q^{|Im(f)|}\]. \\ This enhanced link invariant contains strictly more information than the quandle coloring invariant defined in Definition \ref{countinginva}. \\
    \end{definition}

     In this paper, we will refer to this simply as \textit{quandle coloring polynomial}. We will say that a quandle $(Q,*)$ is able (resp. unable) to distinguish two links $L_1,L_2$ if their quandle coloring polynomials are different (resp. same).\\ 

    \begin{example}
    
The right-handed trefoil knot, denoted by $\overline 3_1$, is illustrated in Figure \ref{trefoil}. 
 \begin{figure}[h]
     \centering
    \includegraphics[width=0.3\linewidth]{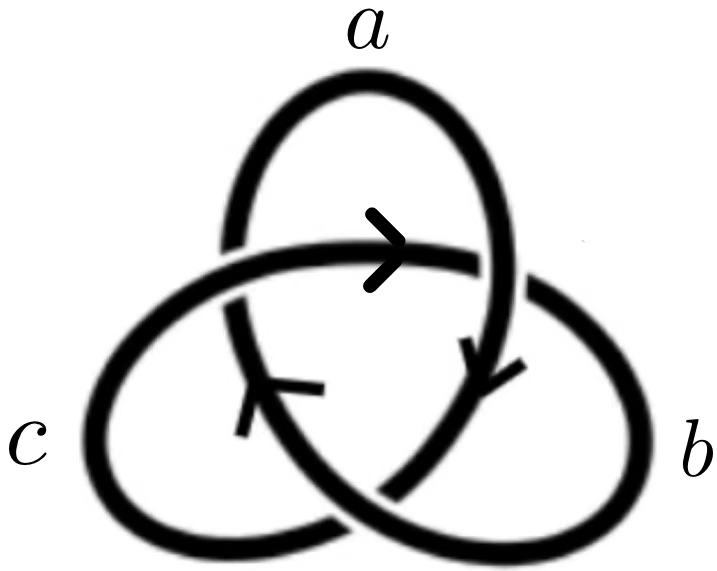}
     \caption{The right-handed trefoil knot.}
     \label{trefoil}
 \end{figure}

The fundamental quandle of $\overline 3_1$ is $Q(\overline 3_1,\times)=\{a,b,c~:~c\times b=a,b\times a=c,a\times c=b\}$. Consider the quandle $\mathbb Z_3( *)$, where $x*y=2x-y,$ for all $x, y \in \mathbb Z_3$. Then deriving the coloring space of this quandle on the trefoil knot is essentially the same as solving this system of equations modulo $3$: $\{2c-b=a,2b-a=c,2a-c=b~(mod~3)\}$. One can check that the solution requires  $a+b+c\equiv 0~(mod~3)$. This means $Hom(Q(\overline 3_1,\times),\mathbb{Z}_3(*))=\{(a\to 1,b\to 1,c\to 1),(a\to 2,b\to 2,c\to 2),(a\to 0,b\to 0,c\to 0), (\sigma(a)\to 0, \sigma(b)\to 1, \sigma(c)\to 2~|~\sigma\in S_{abc}\}$. The coloring counting invariant from $\mathbb{Z}_3(*)$ on the trefoil knot is therefore $|~Hom(Q(\overline 3_1,\times),\mathbb{Z}_3(*))~|=9$. The coloring polynomial is $\Phi(\overline 3_1,\mathbb{Z}_3(*))=3q^1+6q^3$. \\
\end{example}

\section{Medial Quandles and Detecting Causality} \label{section 3}

\begin{figure}[H] 
    \centering
        \includegraphics[width=0.7\linewidth]{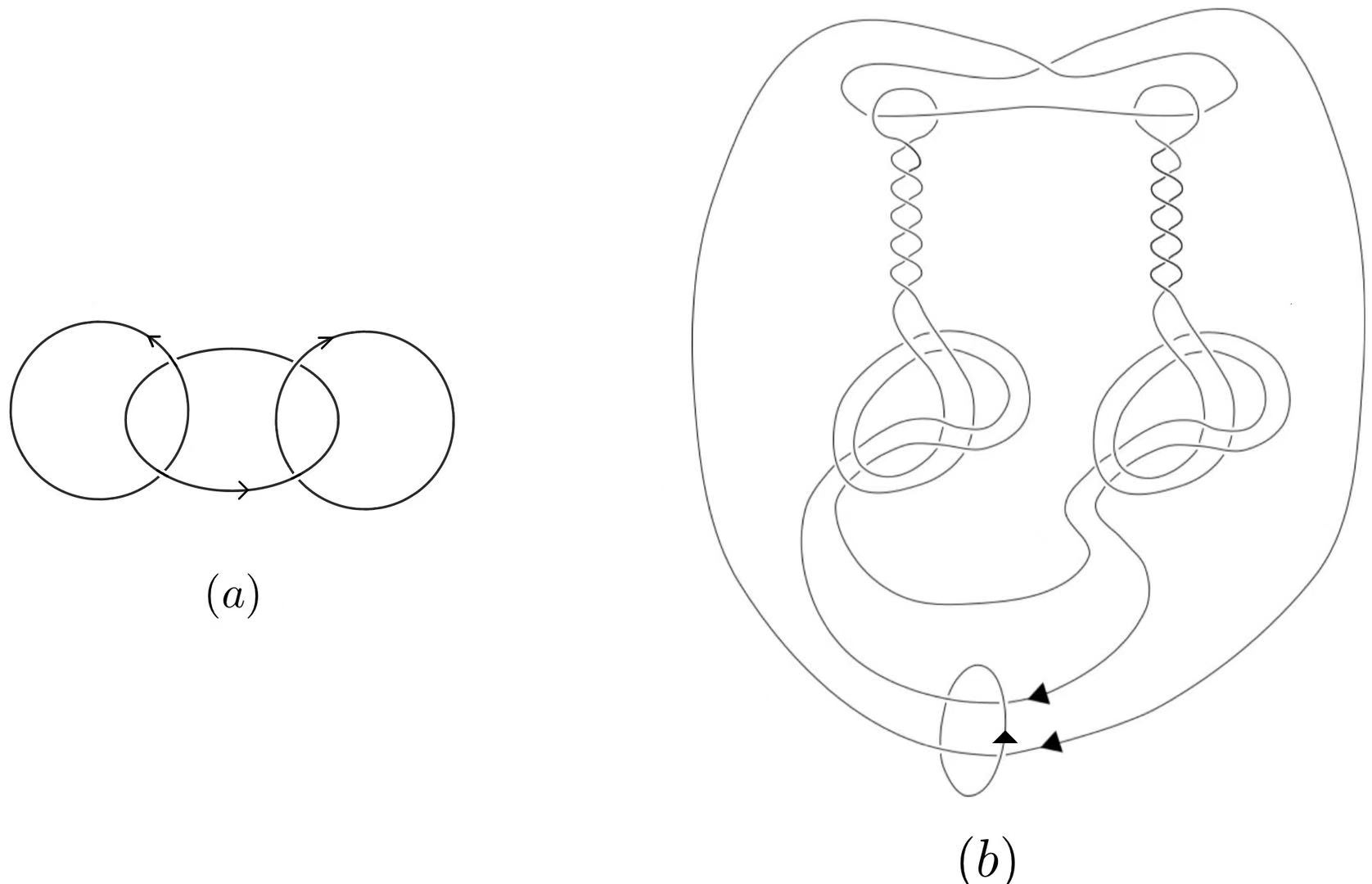}
    \caption{(a): the connected sum of two Hopf links (b): the first Allen-Swenberg link in series, which we denote by $A_1$.} \label{L2HandA1}
\end{figure}

 Let $X$ be a $(2+1)$-dimensional globally hyperbolic spacetime with a Cauchy surface  whose universal cover is homeomorphic to $\mathbb R^2$. Allen and Swenberg \cite{AS} investigated whether there are  link invariants that can detect causality in $X$. They constructed an infinite series of three-component links, which we called Allen-Swenberg links, and suggested that such link invariants must be  able to distinguish between these links and the connected sum of two Hopf links, which we denote by $L_{2H}$ (see Figure \ref{L2HandA1} (a)). They found that the Conway polynomial fails to do so. In this section, we analyze the ability of medial quandles to  distinguish between $L_{2H}$ and the Allen-Swenberg series of links and present main results.

\begin{figure}[H]
  \centering 

      \includegraphics[width=0.82\linewidth]{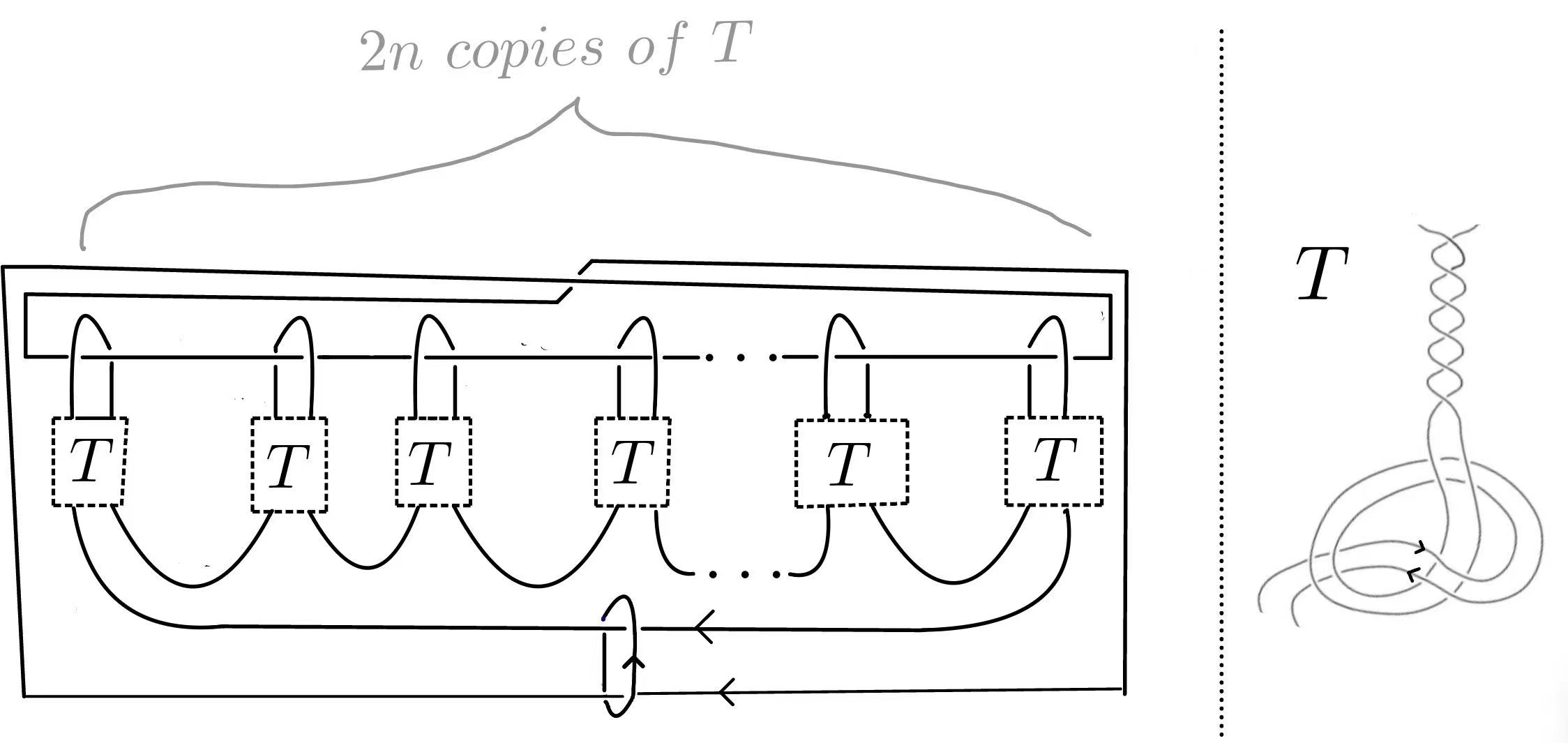}

    \caption{General Allen-Swenberg links $A_n$ and the tangle $T$ that consists of it. } \label{figureAn}\end{figure}

 For a positive integer $n$, the $n^{th}$ Allen-Swenberg link, denoted by $A_n$, is formed by duplicating the tangle in the middle component of $A_1$ as illustrated in Figure \ref{figureAn}.

 We have the following four propositions about the ability of a medial quandle to distinguish $L_{2H}$ and $A_n$. We provide the proofs in Section \ref{sec4}.
\begin{prop} \label{pro 3.1}
    If a medial quandle can distinguish $L_{2H}$ from $A_1$, then it can distinguish $L_{2H}$ from $A_n$, for all $n$. If it cannot distinguish $L_{2H}$ from $A_1$, then it cannot distinguish $L_{2H}$ from $A_n$, for all $n$. \\ \end{prop}  
    \begin{remark}
        
  Proposition \ref{pro 3.1} allows us to treat the infinite series of links $\{A_n\}$ as a single object in the discussion of whether medial quandles can distinguish them from $L_{2H}$. \\ \end{remark}

\begin{prop} \label{prop 3.2}
     Any medial quandle's coloring counting invariant (see Definition \ref{countinginva}) fails to distinguish $L_{2H}$ from $A_n$. Namely, $|Hom(L_{2H},Q(*))|=|Hom(A_n,Q(*))|$. \\ \end{prop}

\begin{prop}\label{prop 3.3}
If a medial quandle $(Q,*)$ can distinguish $L_{2H}$ and $A_n$, then there exists $a,b,c\in Q$ (possibly not distinct) such that 
$a*c=a,~b*c=b,~c*b=c*^{-1}a,$ and $ b*a\neq b$.\\

\end{prop}

We have the following corollary from Proposition \ref{prop 3.3}
\begin{cor}   \label{cor 3}

A medial quandle $(Q,*)$ is unable to distinguish $L_{2H}$ and $A_n$ if the relation $\sim$ on $Q$ defined by $a\sim b\Leftrightarrow a*b=a$ is an equivalence relation. \end{cor} \begin{proof}
Suppose for the sake of contradiction that $Q$ can distinguish between $L_{2H}$ and $A_n$.  Let $a\sim b\Leftrightarrow a*b=a$ be an equivalence equation. Then $\forall a,b,c\in Q$, if $a*c=a$ and $b*c=b$, then $b\sim a$. Therefore, $b*a=b$, which contradicts the inequality in Proposition \ref{prop 3.3}. \\ \end{proof}

The following follows from Corollary \ref{cor 3}:
\begin{cor}
    
 \label{coro_failofAlexanderquandle} Alexander quandles (recall from Example \ref{DofAlexan})  are unable to distinguish $L_{2H}$ from $A_n$. This is because $a*b=a$ implies $(1-f)(x)+f(y)=x\Leftrightarrow f(x)=f(y)$, which is an equivalence relation. \end{cor}

\begin{example} All commutative medial quandles are unable to distinguish $L_{2H}$ from $A_n$. This is because $a*b=a$ implies $b*a=a\to b=a$, which is an equivalence relation. \end{example}

\begin{prop} \label{prop 3.5}
 A medial quandle $(Q,*)$ can distinguish between $L_{2H}$ and $A_n$ if there exists $a,b,c\in Q$ such that\\
$$a*c=a,\ b*c=b,\ c*b=c*^{-1}a,\ \text{and} \ b*(a*^{-1}b)\notin \{a,b,c,c*b\}.$$ \end{prop}

\begin{example}
We construct a $f_{i,j}$ quandle (recall Definition \ref{def_fijquandle}). Let $X_1=\{1,2,4\},~X_2=\{3\}$ and $f_{1,1}=f_{2,2}=f_{2,1}=Id,~f_{1,2}=(1,2)$, where $f_{i,j}:X_i\to X_i$. Let $Q_4=X_1\cup X_2.$ Then the medial quandle $(Q_4,*)$ has the operation table given in Figure \ref{Q4label}. Let $a=3,b=1$, and $c=4$ in Proposition \ref{prop 3.5}. It is easy to check that there does not exist $f\in Hom({L_{2H}},Q_4)$ with $|Im(f)|=4$. However, Figure \ref{Q4label} shows a coloring on $A_1$ by $Q_4$, given by $g\in Hom(A_1,Q_4)$, with $|Im(g|=4$.

\begin{figure}[h] \centering
   \includegraphics[width=0.18\linewidth]{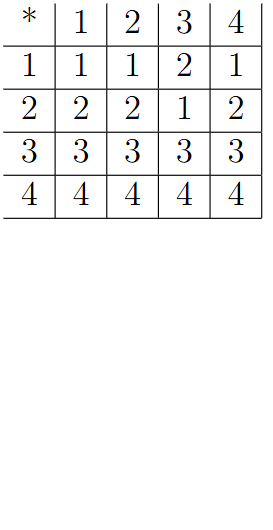}  \includegraphics[width=0.45\linewidth]{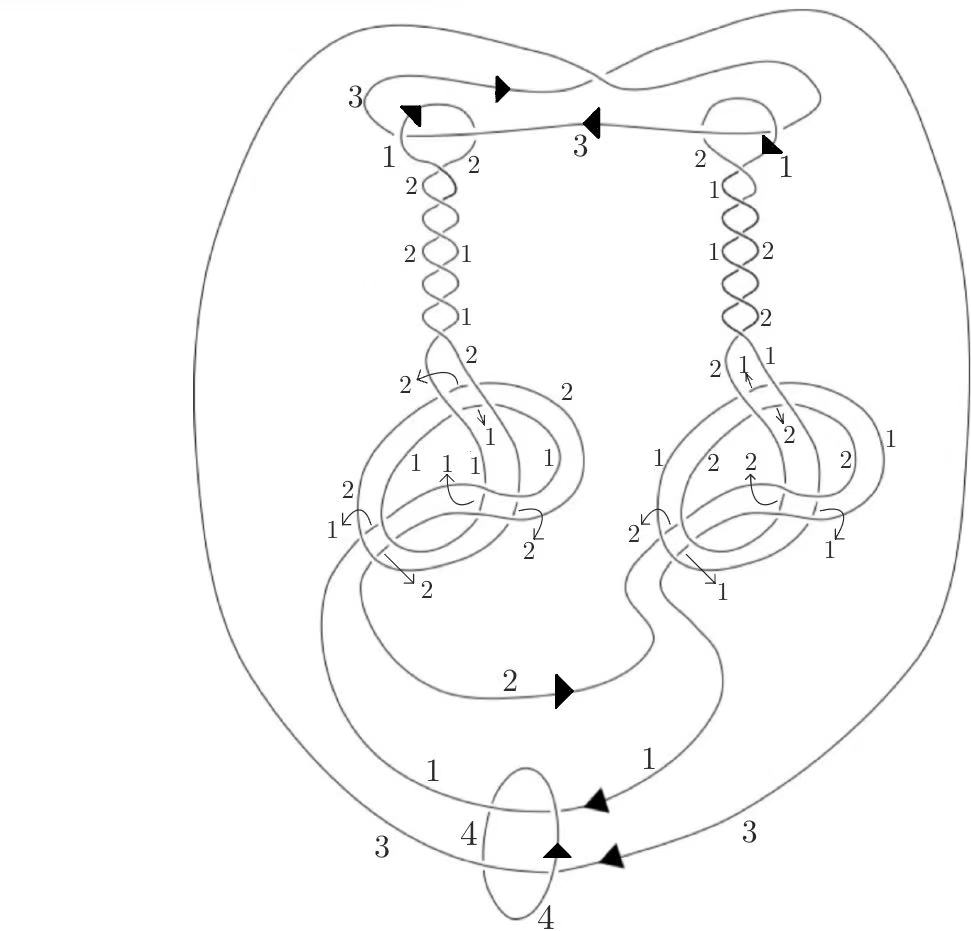} 
   \caption{Operation table of $Q_4$ (left) and the representation of a coloring on $A_1$ by $Q_4$ with 4 elements (right).} \label{Q4label}
\end{figure}
This means that the coloring polynomial from $Q_4$ can distinguish  between $A_1$ and $L_{2H}$ because the polynomial in $A_1$ contains a fourth power term, but the polynomial in $L_{2H}$ does not.
\end{example}

We then discuss a method for constructing quandles that satisfy the condition of Proposition \ref{prop 3.5}.

\begin{prop}
Any quandle obtained by Definition \ref{def_fijquandle} is  medial.
\end{prop}
\begin{proof}
Let $i,j,k,l\in \Lambda$ be arbitrary. Then
\begin{equation} (x_i*x_j)*(x_k*x_l)=f_{i,j}(x_i)*f_{k,l}(x_k)=f_{i,k}(f_{i,j}(x_i))=f_{i,j}(f_{i,k}(x_i))=f_{i,k}(x_i)*f_{j,l}(x_j)=(x_i*x_k)*(x_j*x_l). 
\end{equation} 
\end{proof}

\begin{definition}
A  quandle  is called trivial if $x*y=x$ for all $x,y\in Q$.
\end{definition}
\begin{theorem} \label{coro_successfulfamily}
Let $Q=\{X_i\}_{i\in \Lambda}$ be a non-trivial $f_{ij}$ quandle. Then there exists $Q'$ a medial quandle such that $Q'$ distinguishes  $L_{2H}$ from $A_n$, and $Q$ is a subquandle of $Q'$.
\end{theorem}

\begin{proof}
Since $Q$ is non-trivial, there exists $A,B\in \Lambda,A\neq B$ such that $f_{BA}:X_B\to X_B\neq Id$. We then take arbitrary $b\in X_B$ such that $f_{BA}(b)\neq b$. Take arbitrary $a\in X_A$. Choose an index $C\notin \Lambda$, let $X_C=\{c\}$ and $Q'=Q\sqcup X_C$. Define $f_{iC}=id_{X_i},f_{Ci}=id_{X_C}$, for all $i\in \Lambda$. This constructs a $f_{ij}$ quandle structure on $Q'$ and one can check that the elements $a,b,c\in Q'$ satisfy the equations in Proposition \ref{prop 3.5}.
\end{proof}

 \section{Proofs for Propositions in Section \ref{section 3}} \label{sec4}
 We start by discussing some properties of medial quandle.
 \begin{prop} \label{medial anxioms}
    
 Any medial quandle $Q(*)$, with $*^{-1}$ denoting its inverse operation, satisfies the following nine properties (proofs are included), $\forall x,y,z,x_1,x_2,y_1,y_2\in Q:$\par \begin{enumerate}
 
\item $x*(y*z)=(x*x)*(y*z)=(x*y)*(x*z)$, \item $[(x_1*^{-1}y_1)*(x_2*^{-1}y_2)]*(y_1*y_2)=[(x_1*^{-1}y_1)*y_1]*[(x_2*^{-1}y_2)*y_2)]=x_1*x_2\Rightarrow (x_1*x_2)*^{-1}(y_1*y_2)=(x_1*^{-1}y_1)*(x_2*^{-1}y_2)$,
\item 
Part (2) implies $x*^{-1}(y*z)=(x*x)*^{-1}(y*z)=(x*^{-1}y)*(x*^{-1}z)$,
\item 
Part (2) implies $[(x_1*^{-1}y_1)*^{-1}(x_2*^{-1}y_2)]*(y_1*^{-1}y_2)=x_1*^{-1}x_2$, which in turn implies $ (x_1*^{-1}y_1)*^{-1}(x_2*^{-1}y_2)=(x_1*^{-1}x_2)*^{-1}(y_1*^{-1}y_2)$, \item 
$(x*^{-1}y)*^{-1}z=(x*^{-1}y)*^{-1}(z*^{-1}z)=(x*^{-1}z)*^{-1}(y*^{-1}z)$, \item $x*^{-1}(y*^{-1}z)=(x*^{-1}x)*^{-1}(y*^{-1}z)=(x*^{-1}y)*^{-1}(x*^{-1}z)$, \item 
$(x*^{-1}y)*z=(x*^{-1}y)*(z*^{-1}z)=(x*z)*^{-1}(y*z)$,\item 
$x*(y*^{-1}z)=(x*^{-1}x)*(y*^{-1}z)=(x*y)*^{-1}(x*z)$, and
\item 
$[(x*^{-1}z)*(y*^{-1}z)]*z=x*y\Rightarrow$ $(x*y)*^{-1}z=(x*^{-1}z)*(y*^{-1}z)$  (this is true for all quandles). \end{enumerate} \end{prop}

With the help of Proposition\ref{medial anxioms}, we then prove some further properties of medial quandles  which will be used in the proofs of Propositions \ref{pro 3.1}, \ref{prop 3.2}, \ref{prop 3.3}, and \ref{prop 3.5} later in this section. 

\begin{definition} \label{Dfxy} For any $x,y\in Q$, define $f_{xy}:Q\to Q$ by $f_{xy}(q)=(q*^{-1}x)*y$. Let $f_{x_2y_2}\circ f_{x_1y_1}$ denote the composition of two such functions. \end{definition} 

\begin{prop} \label{propertyfxy}\mbox{} 
Let $(Q,*)$ be a medial quandle. Then

\begin{enumerate}
    \item \label{a1} $[(a*x)*^{-1}y]*z=[(a*z)*^{-1}y]*x,~\forall a,x,y,z\in Q$.\\
\item \label{a2} $[(a*^{-1}x)*y]*^{-1}z=[(a*^{-1}z)*y]*^{-1}x,~\forall a,x,y,z\in Q$.\\ 
\item  \label{a3} $f_{xy}\circ f_{ab}=f_{ab}\circ f_{xy}$. 
    \item \label{a4} $f_{xy}(a*b)=f_{xy}(a)*f_{xy}(b)$ and $f_{xy}(a*^{-1}b)=f_{xy}(a)*^{-1}f_{xy}(b)$. 
    \item  \label{a5}$f_1\circ f_2\circ...\circ f_k=f_{1'}\circ f_{2'}\circ...\circ f_{k'}$, where $\{i\}\to\{i'\}$ is a permutation, and $f_i=f_{x_iy_i}$ for some $x_i,y_i\in Q$.

\end{enumerate} \end{prop}
\begin{proof} \mbox{} \begin{enumerate} 
    \item $[(a*x)*^{-1}y]*z=[(a*x)*z]*^{-1}(y*z)=[(a*z)*(x*z)]*^{-1}(y*z)=[(a*z)*^{-1}y]*[(x*z)*^{-1}z]=[(a*z)*^{-1}y]*x$.
    \item Note that item (4) in Proposition \ref{medial anxioms} implies that $(Q,*^{-1})$ is also a medial quandle. Since \ref{a1} is true, \ref{a2} should be true by symmetry. \item $f_{xy}\circ f_{ab}(q)=[[(q*^{-1}a)*b]*^{-1}x]*y=_1[[(q*^{-1}x)*b]*^{-1}a]*y=_2[[(q*^{-1}x)*y]*^{-1}a]*b=f_{ab}\circ f_{xy}(q)$.
    \item  $f_{xy}(a*b)=[(a*b)*^{-1}x]*y=[(a*^{-1}x)*(b*^{-1}x)]*y=[(a*^{-1}x)*y]*[(b*^{-1}x)*y]=f_{xy}(a)*f_{xy}(b)$\\ $f_{xy}(a*^{-1}b)=[(a*^{-1}b)*^{-1}x]*y=[(a*^{-1}x)*^{-1}(b*^{-1}x)]*y=[(a*^{-1}x)*y]*^{-1}[(b*^{-1}x)*y]=f_{xy}(a)*^{-1}f_{xy}(b)$.
    \item This follows from \ref{a3}.
\end{enumerate}
    
\end{proof}

\begin{theorem}
    
 \label{corofmapsareautomorphism}
The group generated by $\{f_{xy}~|~x,y\in Q\}$ is abelian and its elements are automorphisms on $Q$.
\end{theorem}

\begin{proof}
This follows directly from Proposition \ref{propertyfxy}.
\end{proof}

The following  two lemmas will be used in the proofs of Propositions \ref{pro 3.1}, \ref{prop 3.2}, \ref{prop 3.3}, and \ref{prop 3.5}.
\begin{lemma} \label{lemma1}

Consider a tangle given by a sub-diagram of $A_1$ in Fig \ref{TA}, which is denoted by $T_A$.
       \begin{figure}[h]  \centering       \includegraphics[width=0.57\linewidth]{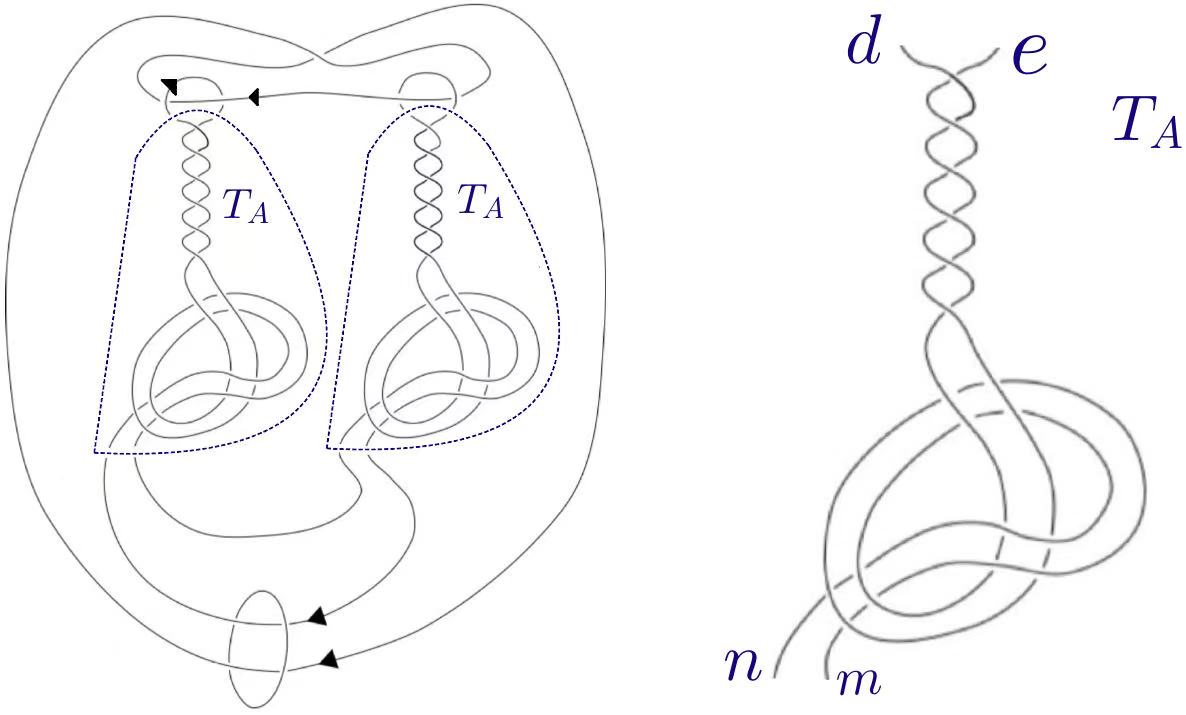} \caption{The tangle $T_A$.} \label{TA}
          \end{figure}
\\For any coloring by medial quandle $(Q,*)$ on $T_A$, where $n,m,d,e,\in Q$ denote elements that color the two pairs of arcs of $T_A$ as illustrated in Fig \ref{TA}, we must have $n=d,m=e$. \end{lemma}

\begin{proof} \label{proflemma1}
     Note that from Figure \ref{TAtransf},  the following tangle is an equivalent form of $T_A$ by applying the third Reidemeister moves  : 
\begin{figure}[h]   \centering  \includegraphics[width=0.73\linewidth]{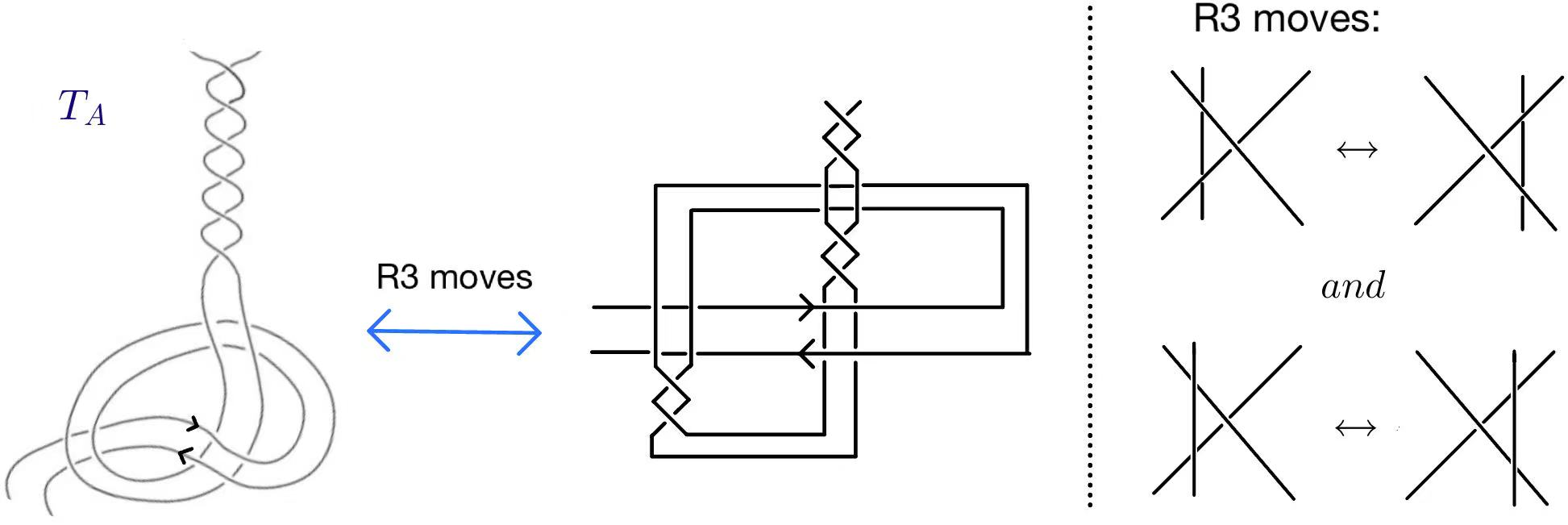} 
\caption{{Alternative diagram of $T_A$.} } \label{TAtransf}
\end{figure}
  
Let $(Q,*)$ be an arbitrary medial quandle. To analyze the coloring of $Q$ on $T_A$ We decompose $T_A$ into two tangles  We first analyze the coloring of $(Q,*)$ on two tangles shown in Figure \ref{fig:colortwostr}, called  $T_1$ and $T_2$, respectively. 
\begin{figure}[H]    \centering
 \includegraphics[width=0.88\linewidth]{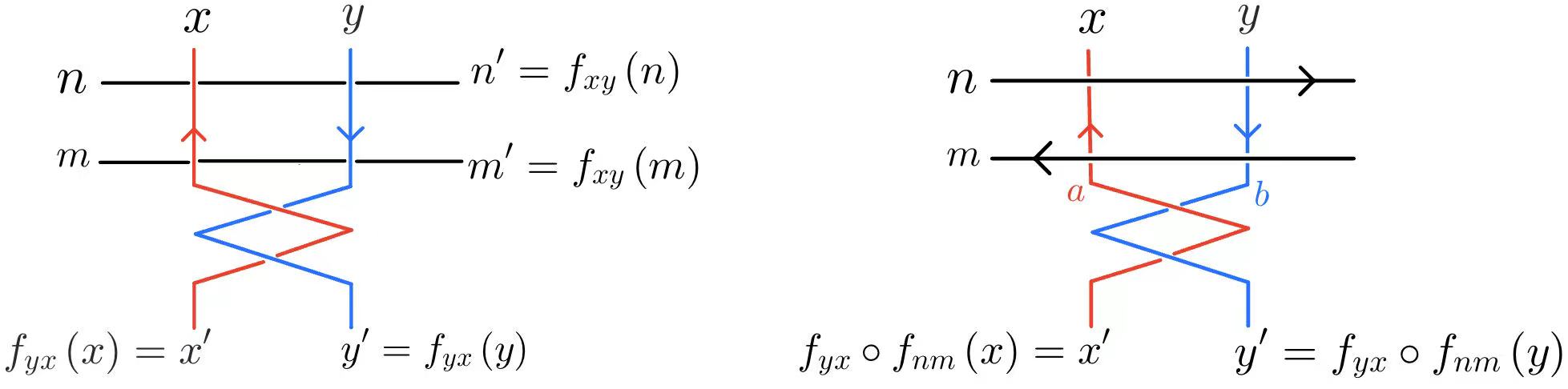} 
 \caption{$T_1$(left) and $T_2$ (right). Recall the definition of $f_{xy}$ in Definition \ref{Dfxy}. }
    \label{fig:colortwostr} \end{figure}
We use letters denoting some elements in $Q$ to label arcs in Fig \ref{fig:colortwostr} to denote a quandle coloring on the structure. For both tangles, the elements that color the four arcs labeled by $x',y',n',m'$, can be expressed by $x,y,n,m$ by applying crossing relation :

For structure \textcircled{1}:\begin{center}
    $n'=(n*^{-1}x)*y=f_{xy}(n)~~m'=(m*^{-1}x)*y=f_{xy}(m)$\\  $y'=y*x=(y*^{-1}y)*x=f_{yx}(y)$\\ $x'=x*^{-1}y'=x*^{-1}(y*x)=(x*^{-1}y)*(x*^{-1}x)=(x*^{-1}y)*x=f_{yx}(x)$ \end{center}

For structure \textcircled{2}: \begin{center}
     $a=f_{nm}(x)~~b=f_{nm}(y)$, so\\
$y'=b*a=f_{nm}(y)*f_{nm}(x)=f_{nm}(y*x)=f_{nm}((y*^{-1}y)*x)=f_{nm}\circ f_{yx}(y)$\\
$x'=a*^{-1}y'=f_{nm}(x)*^{-1}(f_{nm}\circ f_{yx}(y))=f_{nm}(x*^{-1}(y*x))=f_{nm}((x*^{-1}y)*x))=f_{nm}\circ f_{yx}(x)$ \end{center}

 Note that the diagram of $T_A$ in Fig \ref{TAtransf} can be viewed as the combination of $T_1$ and $T_2$:
  \begin{figure}[H] \centering \includegraphics[width=0.45\linewidth]{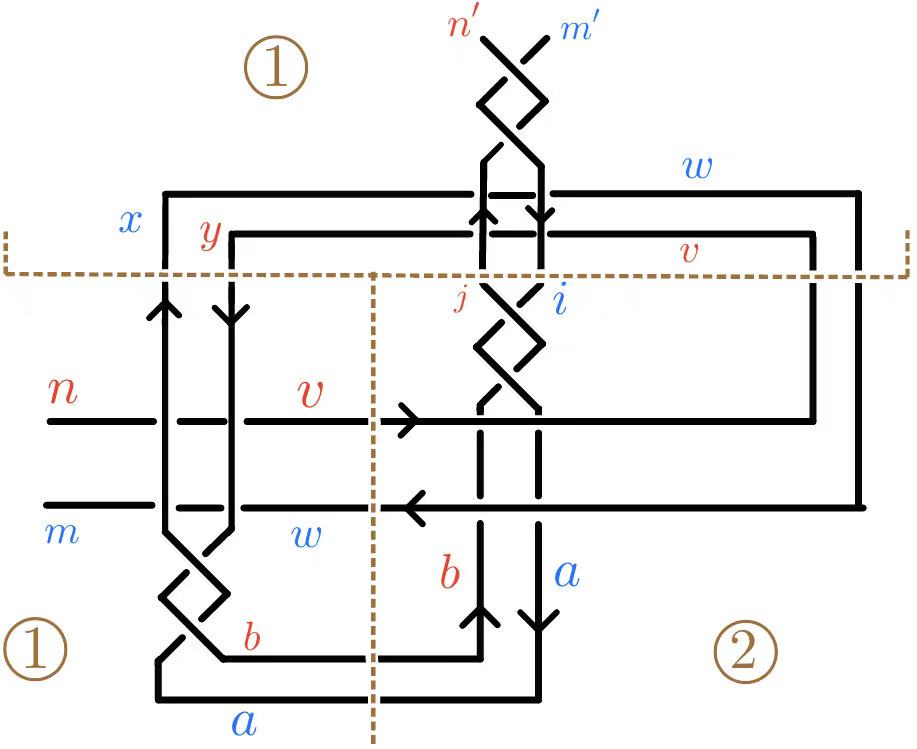} 
  \caption{Decomposition of $T_A$.} \label{TAascombine} \end{figure}

  By applying the calculations from Figure \ref{fig:colortwostr} on the labels in Figure \ref{TAascombine}, we have \begin{equation}
       n'=f_{ji}(j),~j=f_{wv}\circ f_{ba}(b),~b=f_{yx}(y),~y=f_{ij}(v),~v=f_{xy}(n), \text{and} \label{lemma4.3 eq1}
  \end{equation}
  \begin{equation}
      m'=f_{ji}(i),~i=f_{wv}\circ f_{ba}(a),~a=f_{yx}(x),~x=f_{ij}(w),~w=f_{xy}(m). \label{lemma4.3 eq2}
  \end{equation}
   \\To simplify the expression for $j$ and $i$, we derived from equations \ref{lemma4.3 eq1} and \ref{lemma4.3 eq2} that 
  \begin{center}
       
   $b=f_{yx}\circ f_{ij}(v),~a=f_{yx}\circ f_{ij}(w)$.  \end{center} Let $g:=f_{yx}\circ f_{ij}$. By Theorem \ref{corofmapsareautomorphism}, $g$ is an automorphism on $Q$ and commutes with all $f_{xy}$-type maps. Therefore, 
   \begin{equation}
   j=f_{wv}\circ f_{ba}(b)=f_{wv}[(g(v)*^{-1}g(v))*g(w)]=f_{wv}(g(v*w))=g(f_{wv}(v*w))=g(v)=b. \end{equation} 
   \begin{equation} i=f_{wv}\circ f_{ba}(a)=f_{wv}[(g(w)*^{-1}g(v))*g(w)]=f_{wv}(g((w*^{-1}v)*w))=g(f_{wv}((w*^{-1}v)*w))=g(w)=a. \end{equation}
   Plugging in $j=b,~i=a$ to equations  \ref{lemma4.3 eq1} and \ref{lemma4.3 eq2} gives \begin{equation}
   n'=f_{ji}\circ f_{yx}\circ f_{ij}\circ f_{xy}(n), m'=f_{ji}\circ f_{yx}\circ f_{ij}\circ f_{xy}(m). \end{equation} Note that $f_{pq}\circ f_{qp}$ is the identity function for all $p,q\in Q(*)$ (recall Definition \ref{Dfxy}). By rearranging each individual functions which item \ref{a5} in Proposition \ref{propertyfxy} allows, we have\\
  $n'=f_{ji}\circ f_{ij}\circ f_{yx}\circ f_{xy}(n)=f_{ji}\circ f_{ij}(n)=n$ and $m'=f_{ji}\circ f_{ij}\circ f_{yx}\circ f_{xy}(m)=f_{ji}\circ f_{ij}(m)=m$, as desired. \end{proof}

\begin{lemma} \label{lemma2}
Given a medial quandle $(Q,*)$ and for any $a_0,b_0\in Q$, there exists a unique coloring on the tangle $T_A$ such that the top two arcs (labeled in Fig \ref{topcolorTA}) are colored by $(a_0,b_0)$. \end{lemma} 
\begin{proof}
    
Referring to Figure \ref{topcolorTA}, note that the coloring of the top two arcs of $T_A$ (labeled by $a_0,b_0$) completely determine the coloring on the upper parts of $T_A$ with the 6 twisted crossing relations on the top, so we can fix $a,b\in Q$ that color the two arcs below the 6th crossing.\\
\begin{figure}[H]
    \centering
    \includegraphics[width=0.21\linewidth]{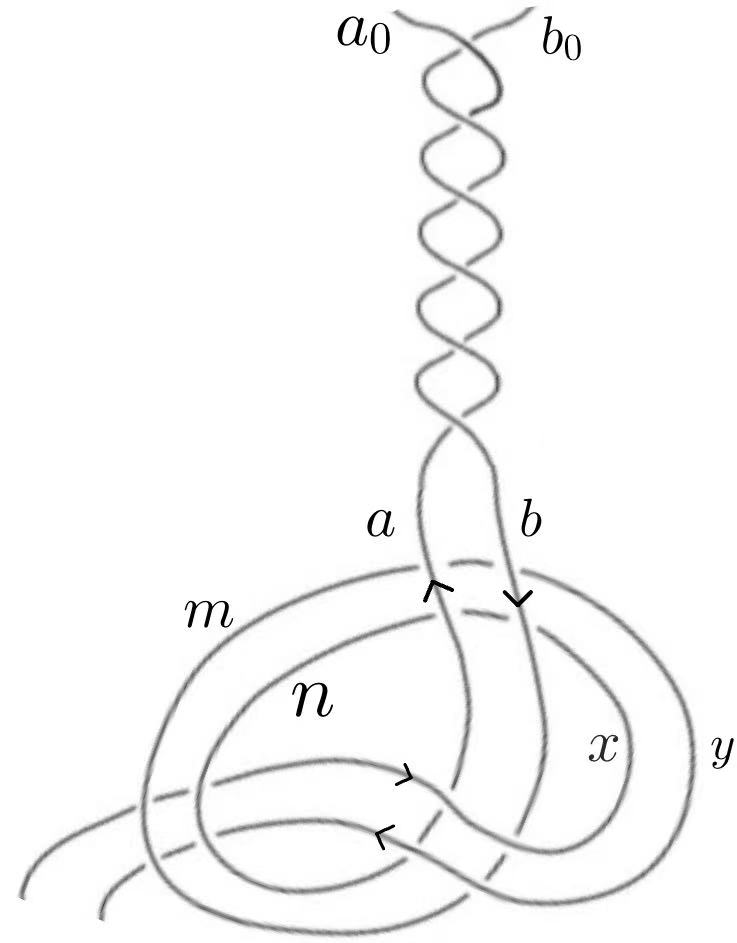} \caption{$T_A$ with arcs partially labeled. } \label{topcolorTA}
\end{figure}

Consider $m,n,x,y$, the coloring elements of four of the arcs of $T_A$, as labeled in Figure \ref{topcolorTA}. These four elements completely determine the coloring on other parts of $T_A$, so to prove Lemma \ref{lemma2} it suffices to show that given fixed $a,b$ there exist unique $m,n,x,y\in Q$ that satisfy all the crossing relations of $T_A$. The crossing relations between them give: \begin{equation}
 (a,b)=(f_{yx}(n),f_{yx}(m))=(f_{yx}\circ f_{ba}(x),f_{yx}\circ f_{ba} (y))=(f_{ba}\circ f_{yx}(x),f_{ba}\circ f_{yx} (y)) \end{equation} 
 \begin{equation} \Leftrightarrow f_{ab}(a)=f_{yx}(x), f_{ab}(b)=f_{yx}(y)\Leftrightarrow  a*b=x*^{-1}(y*x), y=f_{ab}(b)*^{-1}x \end{equation}
 \begin{equation} \Leftrightarrow (a*b)*f_{ab}(b)=x, f_{ab}(b)*^{-1}[(a*b)*f_{ab}(b)]=y. \end{equation}
Hence, $a,b$ completely determine $x,y$, and also $m,n$ by $(n,m)=(f_{ba}(x),f_{ba}(y))$, which completes the proof. \end{proof}

\subsection{Condition on Coloring on $A_1$ and $A_n$} \label{equivalent}
Recall that the $n^{th}$ Allen-Swenberg link is denoted by $A_n$. An arbitrary coloring on $A_1$ by a medial quandle $(Q,*)$ can be represented by Figure \ref{equiA1}:

  \begin{figure}[h] \centering \includegraphics[width=0.5\linewidth]{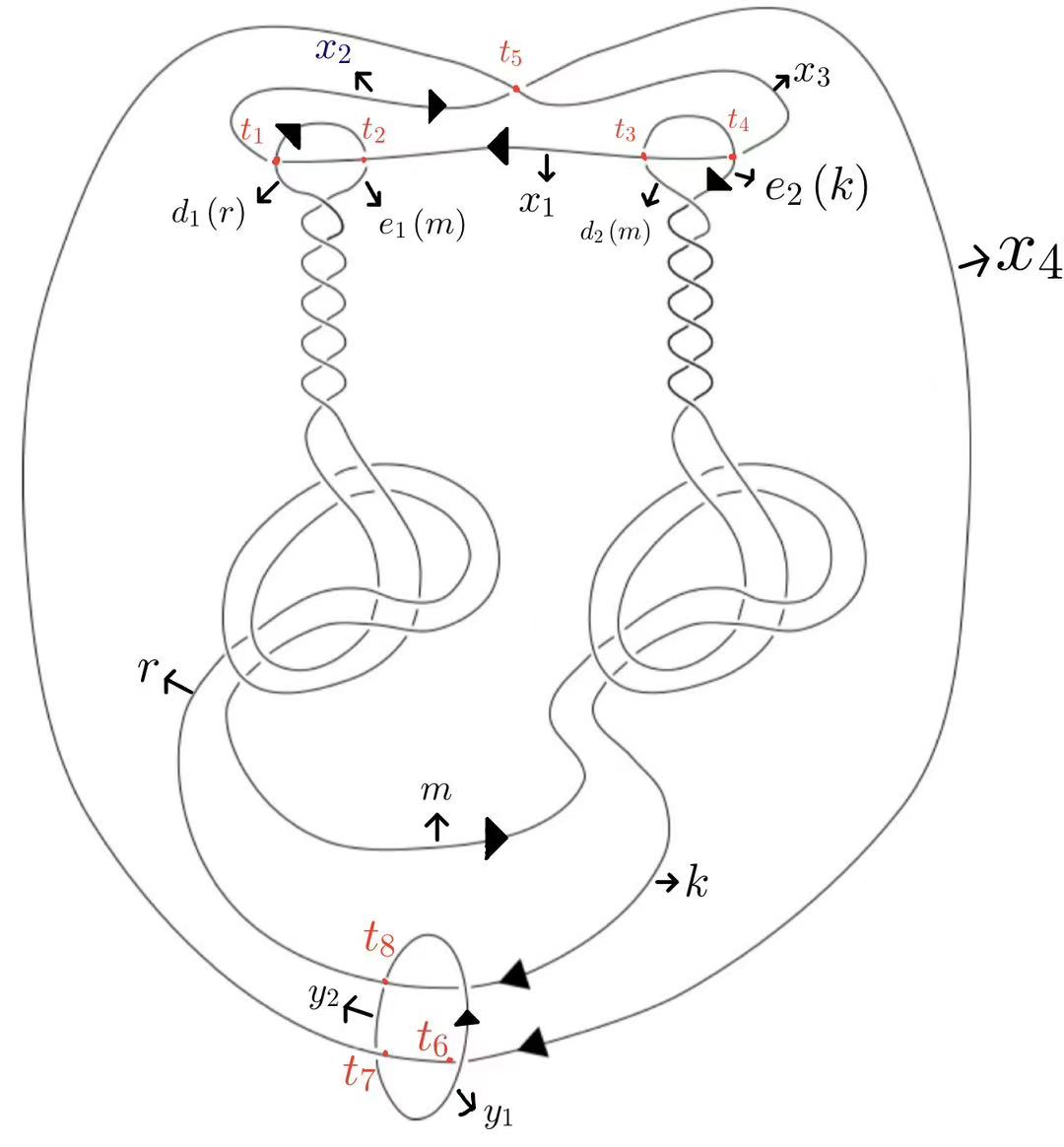}
    \caption{$A_1$ with labeling on some arcs to represent a coloring. Some crossings are also labeled.} \label{equiA1} \end{figure}
We apply Lemma \ref{lemma1} to derive $d_1=r,e_1=d_2=m,$ and $e_2=k.$ The crossing relations that come from the top part of $A_1$ are:
    
\begin{equation} {t_2~and~t_3}: r=m*^{-1}x_1=k,~~~{t_1~and~t_4}: x_2=x_1*r=x_1*k=x_3,~~~{t_5}: x_4=x_2*x_3=x_3*x_3=x_3. \end{equation} Each {$t_i$} refers to a crossing as labeled in Figure \ref{equiA1}, and the equations at right hand side are the crossing relation they give. So in conclusion, \begin{equation}
r=k,x_2=x_3=x_4. \label{eq I} \end{equation} 
By applying Lemma \ref{lemma1} on each $L_i$, an arbitrary coloring on the general Allen-Swenberg link $A_n$ is of the form:
\begin{figure}[H]  \centering   \includegraphics[width=0.75\linewidth]{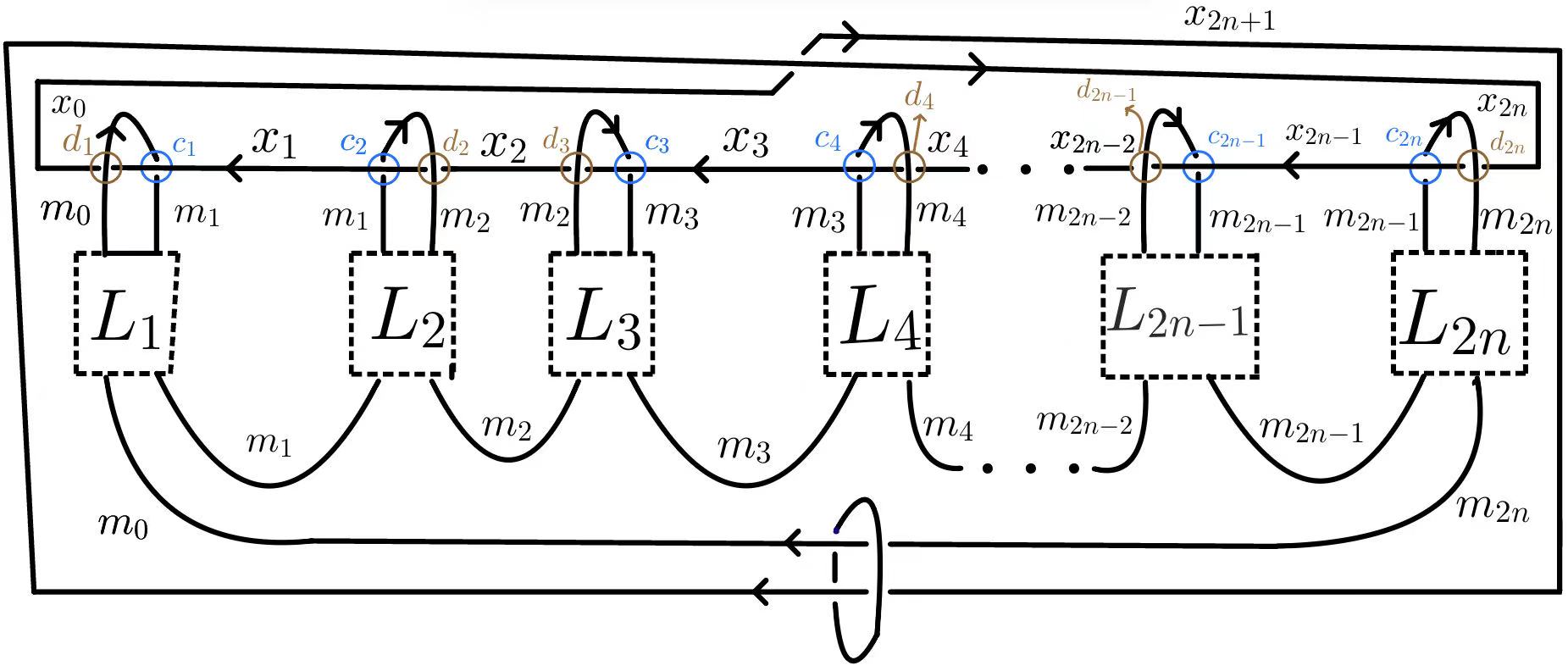} \caption{General coloring on $A_n$ derived by applying Lemma \ref{lemma1} to each $L_i$.}

\end{figure}
For all $i\in \{1,3,...,2n-1\}$, the pair of crossings {$c_i$} and {$c_{i+1}$} require $m_{i-1}=m_i*^{-1}x_i=m_{i+1}$. Therefore, $m_0=m_2,m_2=m_4,...$, so we have  \begin{equation} m_0=m_2=m_4=...=m_{2n}. \label{eq 1} \end{equation} Combining with Equation \ref{eq 1}, the pair of crossings {$d_i$} and {$d_{i+1}$} give \\$x_{i-1}=x_i*m_{i-1}=x_i*m_{i+1}=x_{i+1}$, so similarly \begin{equation} x_0=x_2=x_4=...=x_{2n}. \label{eq 2}\end{equation} For each $j\in\{2,4,...,2n\}$,the crossing {$d_{j}$} gives $x_{j-1}=x_j*^{-1}m_j$. Combining this with Equations \ref{eq 1} and \ref{eq 2}, we have $x_{j+1}$ is the same for all $j$, namely \begin{equation} x_1=x_3=...=x_{2n-1}. \label{eq 3}\end{equation} Each crossing {$c_{j-1}$} requires $m_{j-1}=m_{j-2}*x_{j-1}$ which, combined with Equations \ref{eq 1} and \ref{eq 3}, implies \begin{equation} m_1=m_3=...=m_{2n-1}. \label{eq 4} \end{equation}

\subsection{Proof of Proposition \ref{pro 3.1}}
\begin{figure}[H]
  \centering \includegraphics[width=0.71\linewidth]{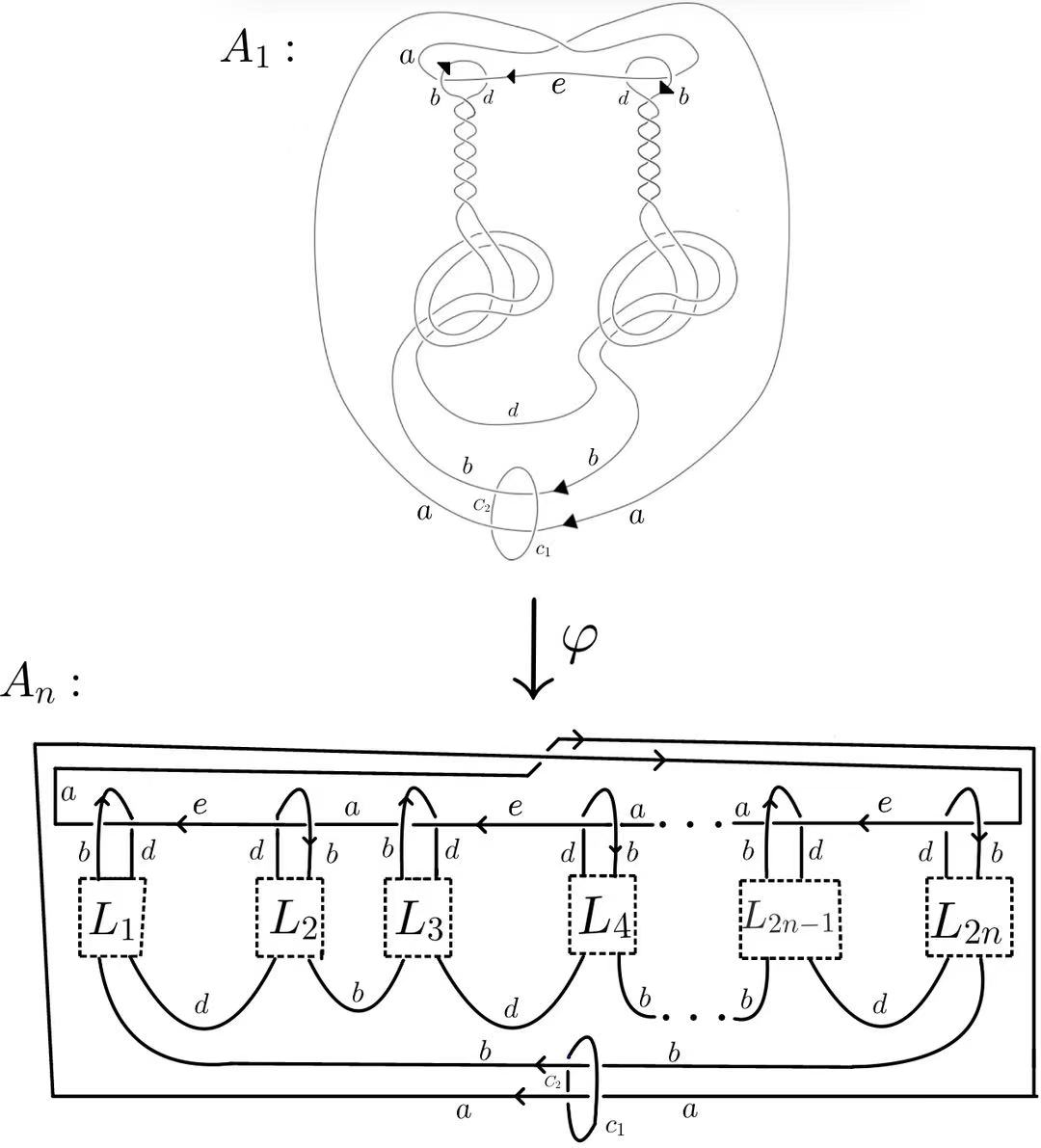} \caption{A bijective map from coloring space on $A_1$ to coloring space on $A_n$ that preserves the number of coloring elements.  Using Lemma \ref{lemma2}, the coloring of unlabeled part of $A_n$ are defined by the unique coloring consistence with the labeling $b$ and $d$. } 
  \label{phi}\end{figure}

To prove Proposition \ref{pro 3.1}, it suffices to show that $\Phi(A_1,Q)=\Phi(A_n,Q)~\forall n\in\mathbb{Z}_+$, for any medial quandle $Q$. We do this by constructing a bijective function $\varphi:Hom(A_1,Q)\to Hom(A_n,Q)$ in Figure \ref{phi}, which satisfies $|Im(h)|=|Im(\varphi(h))|~\forall h\in Hom(A_1,Q)$.

This map is well defined because any coloring on $A_1$ must be of the form shown in Figure \ref{phi} for some $a,b,c_1,c_2,d,e\in Q$  by using Relation \ref{eq I} in Subsection \ref{equivalent}. It is easy to check that $\varphi$ is injective. It is surjective because by Relations \ref{eq 1}, \ref{eq 2}, \ref{eq 3} and \ref{eq 4} in Subsection \ref{equivalent}, any coloring on $A_n$ must be  of the form shown in Figure \ref{phi}. Also, $\forall h\in Hom(A_1,Q), Im(h)=Im(\varphi(h))$ so $|Im(h)|=|Im(\varphi(h))|$. Hence,\\

$$\Phi(A_1,Q)=\sum_{h\in Hom(A_1,Q)} q^{|Im(h)|}=\sum_{g\in Hom(A_n,Q)} q^{|Im(g)|}=\Phi(A_n,Q),$$ as desired.

\subsection{Proof of Proposition \ref{prop 3.2}}
\begin{figure}[H] \centering \includegraphics[width=0.75\linewidth]{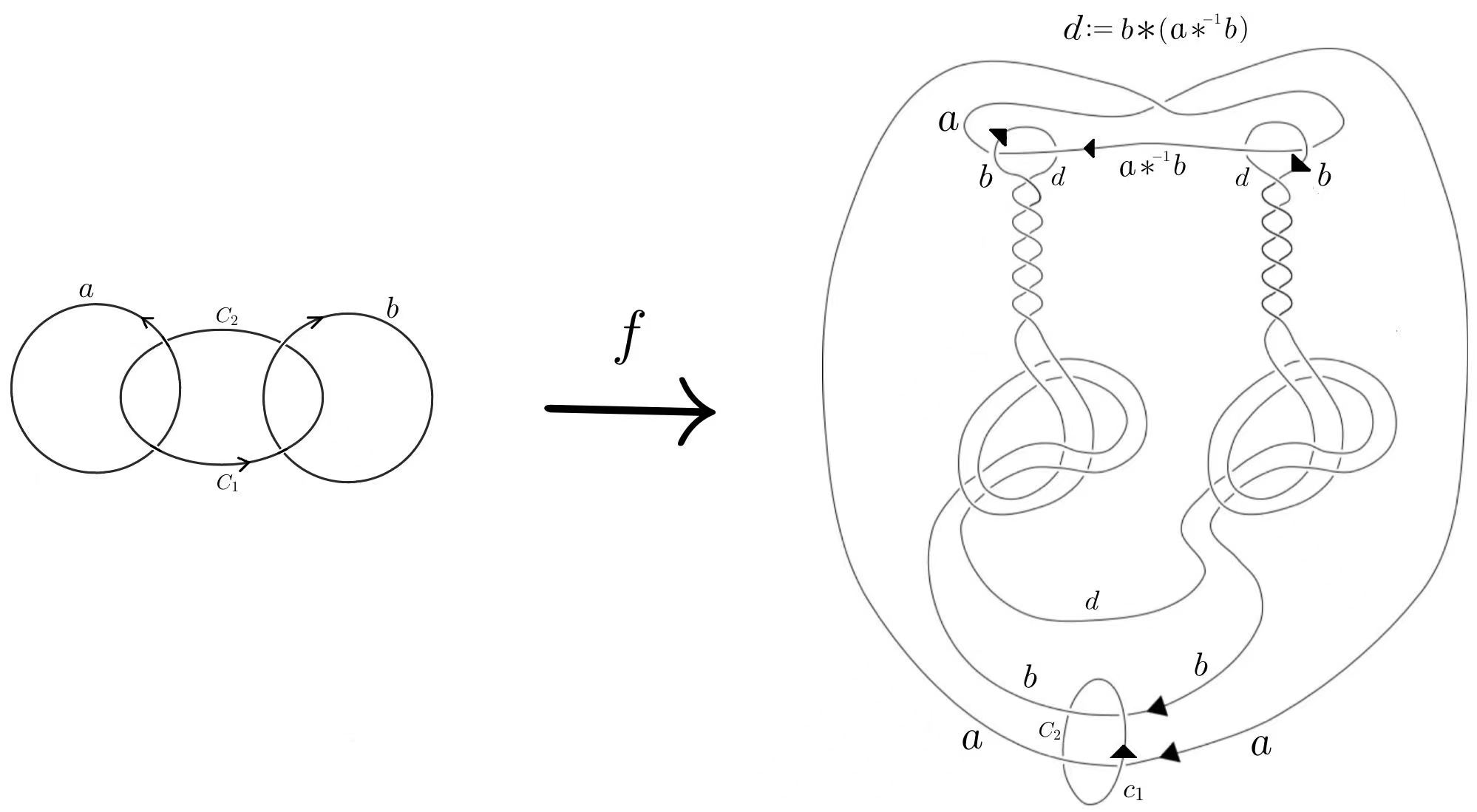}
     \caption{$f:Hom(L_{2H},Q)\to Hom(A_1,Q)$.The unlabeled parts of $A_1$ are colored by the unique coloring determine by $(b,d)$ and  $(d,b)$ respectively using Lemma \ref{lemma2}.}\label{f:L2HtoA1}\end{figure}

We construct a function $f:Hom(L_{2H},Q)\to Hom(A_1,Q)$ as shown in Figure \ref{f:L2HtoA1}.\\

 The map $f$ is obviously injective. It is also surjective  by Relation \ref{eq I} in Subsection \ref{equivalent} and Lemma \ref{lemma1}. Hence, $f$ is bijective, as desired.

\subsection{Proof of Proposition \ref{prop 3.3}}
\begin{figure}[H]  \centering  \includegraphics[width=0.8\linewidth]{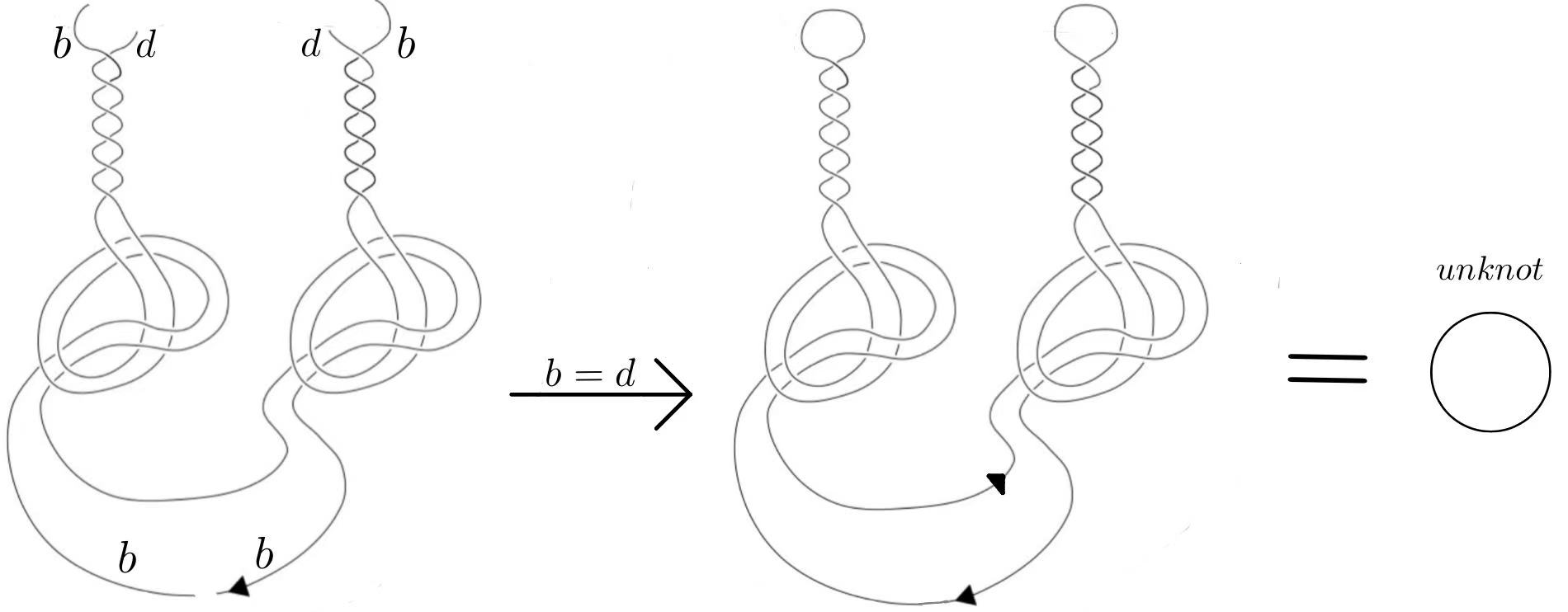} \caption{The middle component of $A_1$ become the same as the unknot when equalizing the two ending pairs of arcs. } \label{middletrans}\end{figure}

By Proposition \ref{pro 3.1}, it suffices to show that a medial quandle $(Q,*)$ is unable to distinguish $L_{2H}$ from $A_1$ if $a*c=a,~b*c=b,~c*b=c*^{-1}a$ implies $b*a=b$ for all $a,b,c\in Q$. \\Note that these three equalities are exactly the crossing relations of $L_{2H}$ in Figure \ref{f:L2HtoA1} (left) (by letting $c=c_1$). Referring to Figure \ref{f:L2HtoA1} (right), $d=b*(a*^{-1}b)=(b*a)*^{-1}(b*b)=b*^{-1}b=b$ for any coloring on $A_1$. When added with the condition $d=b$, crossing relations of the complex middle component of $A_1$ become the same as that of a knot equivalent to the unknot, as shown in Figure \ref{middletrans}.

This means that all coloring elements from a quandle on this component have to be the same, which is $b$. This means $f$ defined in Figure \ref{f:L2HtoA1} is not only bijective, but also satisfies  $|Im(h)|=|Im(f(h))|~\forall h\in Hom(L_{2H},Q)$, which means $\Phi(L_{2H},Q)=\Phi(A_1,Q)$, as desired.

\subsection{Proof of Proposition \ref{prop 3.5}} \label{3.5 Proof}
For simplicity, let $H_L:=Hom(L_{2H},Q) , H_A:= Hom(A_1,Q)$, $ L_n:=\{h\in H_L:|Im(h)|=n\}$ and $k_n:=\{h\in H_A:|Im(h)|=n\}$ for $n\in \mathbb{Z_+}$. Suppose for the sake of contradiction that  $\Phi(L_{2H},Q)=\Phi(A_1,Q)$. Then $|L_n|=|k_n|=0~ \forall n>4$ (because $max\{|Im(h)|:h\in H_L\}=4$) and $|L_n|=|k_n|$ for $n=1,2,3,4$. Note that $f$ defined in Figure \ref{f:L2HtoA1} satisfies $|Im(f(h))|\geq |Im(h)|, \forall h\in H_L$. Since $f$ is injective, we must have $f(L_4)=k_4$, which implies $\{h\in H_L\setminus L_4:f(h)\in k_4\}=\emptyset$. This then implies $f(L_3)=k_3$, so similarly we have $\{h\in H_L\setminus L_3:f(h)\in k_3\}=\emptyset$. By repeating this process, we show that $f(L_n)=k_n ~\forall n=1,2,3,4$. This means that $|Im(f(h))|=|Im(h)|, \forall h\in H_L$. However, by choosing $a,b,c$ to be the elements that satisfy the conditions in the statement of Proposition \ref{prop 3.5} and letting $h\in H_L$ be the coloring represented by the left part of Figure \ref{f:L2HtoA1}, we have $Im(h)=\{a,b,c_1,c_2\}$ and $a,b,c_1,c_2,b*(a*^{-1}b)\in Im(f(h))$, so $|Im(f(h))|\geq |Im(h)|+1$, which is a contradiction. This completes the proof.

\section{Conclusion} \label{secconclusion}
The ability of a link invariant to distinguish the Allen-Swenberg links from the connected sum of two Hopf links implies possible capability for it to detect causality in spacetime, and inability to distinguish them implies inability to detect causality in spacetime. We use  algebraic methods to obtain a comprehensive analysis of the ability of medial quandles to complete this task. Specifically, Propositions \ref{prop 3.3} and \ref{prop 3.5} give rules to determine whether a given medial quandle is capable of distinguishing the links that are close to an if and only if condition. In particular, this rules out the possibility for the Alexander quandles to complete this task (Corollary \ref{coro_failofAlexanderquandle}). We also obtain a method of constructing a large family of medial quandle that are able to distinguish the connected sum of two Hopf links and every link in the Allen-Swenberg infinite series using Theorem \ref{coro_successfulfamily}, which therefore give powerful candidates for detecting causality in spacetime.

\section{Acknowledgements}
This project was originally initiated as follow-up work for the Horizon Academic Research Program in Math in summer 2022. I shall thank my mentors Vladimir Chernov and Emanuele Zappala for setting up the background and motivation for this project.\\
I would like to express particular gratitude to Professor Rhea Palak Bakshi at UC Santa Barbara, who offered continuous and long-term guidance while I was writing this paper. Her help during my sophomore and junior years providing knot theory knowledge, suggestions on formalizing the proofs and language in the paper is important for this output .


\begin{thebibliography}{99}

\bibitem{CN} 
 V. Chernov and S. Nemirovski. Legendrian links, causality, and the low conjecture. {\it Geom. Funct. Anal.}, 19 (0222503): 1323–1333, 2010. \href{https://arxiv.org/abs/0810.5091}{arXiv:0810.5091} [math.SG].
 \bibitem{L} R. J. Low. Causal relations and spaces of null geodesics. PhD thesis, Oxford University, 1988. \label{L}
\bibitem{MP} S. Mukherjee and J.~H. Przytycki, On the rack homology of graphic quandles, in {\it Nonassociative mathematics and its applications}, 183--197, Contemp. Math., 721, Amer. Math. Soc., RI, ; MR3898509. \label{MP}

\bibitem{AS} S. Allen and J.~H. Swenberg, Do link polynomials detect causality in globally hyperbolic spacetimes?, J. Math. Phys. {\bf 62} (2021), no.~3, Paper No. 032503, 16 pp.; MR4227566. \label{AS}
\bibitem{J}D. Joyce, A classifying invariant of knots, the knot quandle, J. Pure Appl. Algebra {\bf 23} (1982), no.~1, 37--65; MR0638121. \label{J}
\bibitem{NN}E.~A. Navas and S. Nelson, On symplectic quandles, Osaka J. Math. {\bf 45} (2008), no.~4, 973--985; MR2493966. \label{NN}
\bibitem{NT} J. Nat\'ario and P. Tod, Linking, Legendrian linking and causality, Proc. London Math. Soc. (3) \textbf{88} (2004), no. 1, 251--272; MR2031422.
       
     \end{thebibliography}
\end{document}